\documentclass{pnastwo}

\usepackage{pstricks, pst-node, amssymb}

\usepackage{rotating}
\newcommand{\rotxc}[1]{\begin{sideways}#1\end{sideways}}
\newcommand{\invert}[1]{\rotxc{\rotxc{#1}}}

\psset{unit=1pt, arrowsize=4pt, linewidth=.7pt}
\psset{linecolor=blue}
\newgray{grayish}{.90}
\newrgbcolor{embgreen}{0 .5 0}
\def\Le{\invert{$\Gamma$}}
\def\vblack(#1, #2)#3{\cnode*[linecolor=black](#1, #2){3}{#3}}
\def\vwhite(#1,#2)#3{\cnode[linecolor=black,fillcolor=white,fillstyle=solid](#1,#2){3}{#3}}
\countdef\x=23
\countdef\y=24
\countdef\z=25
\countdef\t=26

\def\tbox(#1,#2)#3{
\x=#1 \y=#2 
\multiply\x by 12 
\multiply\y by 12 
\z=\x \t=\y
\advance\z by 12 
\advance\t by 12 
\psline(\x,\y)(\x,\t)(\z,\t)(\z,\y)(\x,\y)
\advance\x by 6
\advance\y by 6 
\rput(\x,\y){{\bf #3}}}

\usepackage{amssymb}
\usepackage{graphicx}
\usepackage{epstopdf}
\usepackage{amscd}
\usepackage{graphics}

\let\trueint=\int
\let\truesum=\sum
\def\int{\mathop{\textstyle\trueint}\limits}
\def\sum{\mathop{\textstyle\truesum}\limits}
\let\<=\langle
\let\>=\rangle

\def\S{{\mathcal S}}

\usepackage{amsmath,  amssymb, amsbsy}
\usepackage{amsfonts, latexsym, stmaryrd, amscd, xy}
\usepackage[mathscr]{eucal}
\usepackage{epsfig}
\newtheorem{example}[theorem]{Example}
\newtheorem{algorithm}[theorem]{Algorithm}
\usepackage{amsfonts}
\usepackage{xy}
\usepackage{amssymb}
\usepackage{amsmath}

\newcommand{\R}{\mathbb R}

\newcommand{\Grkn}{(Gr_{kn})_{\geq 0}}
\newcommand{\Grkntop}{(Gr_{kn})_{> 0}}

\DeclareMathOperator{\CC}{\mathcal C}

\DeclareMathOperator{\M}{\mathcal M}
\DeclareMathOperator{\SSS}{\mathcal S}

\newcommand{\thmrefer}[1]{\renewcommand\thetheorem
  {\protect\ref{#1}}\addtocounter{theorem}{-1}}

\xyoption{all}
\CompileMatrices

\begin{document}
\title{KP solitons, total positivity, and cluster algebras}
\author{Yuji Kodama
\affil{1}{Department of Mathematics, Ohio State University, Columbus, OH 43210}
\and
Lauren Williams 
\affil{2}{Department of Mathematics,
University of California, Berkeley, CA 94720}}


\maketitle

\begin{article}

\begin{abstract}
Soliton solutions of the KP equation have been studied since 1970, 
when Kadomtsev and Petviashvili proposed a two-dimensional nonlinear dispersive
wave equation now known as the KP equation.
It is well-known that
the Wronskian approach to the KP equation provides a method to construct soliton solutions.
The {\it regular} soliton solutions that one obtains in this way
come from points of the totally non-negative part of the 
Grassmannian.  In this paper we explain how the theory of total positivity
and cluster algebras provides a framework for understanding these
soliton solutions to the KP equation.
We then use this framework to
give an explicit construction of certain soliton contour graphs, and
solve the inverse problem for soliton
solutions coming from the totally positive part of the Grassmannian.
\end{abstract}

\maketitle


\section{Introduction}

The {\it KP equation}, introduced in 1970 \cite{KP70}, is
considered to be a prototype of an integrable nonlinear 
dispersive wave equation with two spatial dimensions.
Concretely, solutions to this equation
provide a close approximation to the behavior of shallow water waves, such as beach waves.
Given a point $A$ in the real Grassmannian, one can construct a solution
to the KP equation \cite{Sato};
this solution $u_A(x,y,t)$ is written in terms of a $\tau$-{\it function}, which is
a sum of exponentials. 
More recently, several authors \cite{K04, BC06, CK1,CK3} 
have focused on understanding
the {\it regular} soliton solutions that one obtains in this way:
these come from points of the totally non-negative part of the Grassmannian.

The classical theory of total positivity concerns square matrices in which all minors
are positive.  This theory was pioneered in the 1930's by Gantmacher, Krein, and Schoenberg, 
and subsequently generalized in 
the 1990's by Lusztig \cite{Lusztig2, Lusztig3}, who in particular introduced
the totally positive and non-negative parts of real partial flag varieties.  

One of the most important partial flag varieties is the Grassmannian.
Postnikov \cite{Postnikov} investigated the
totally non-negative part of the
Grassmannian $(Gr_{kn})_{\geq 0}$, which can be 
defined as the subset of the real Grassmannian where all Pl\"ucker
coordinates are non-negative.  
Specifying which minors are strictly positive
and which are zero gives a decomposition into {\it positroid cells}.  
Postnikov introduced a variety of combinatorial objects, including
{\it decorated permutations, $\Le$-diagrams, plabic graphs,} and 
{\it Grassmann necklaces,}
in order to index the cells and describe their properties.  

In this paper we develop a tight connection between 
the theory of total positivity for the Grassmannian
and the behavior of the corresponding soliton solutions to the KP equation.  
To understand a soliton solution $u_A(x,y,t)$, one fixes 
the time $t$, and plots the points where $u_A(x,y,t)$ has a local maximum.  
This gives rise to a {\it tropical curve} in the $xy$-plane; concretely, this
shows the positions in the plane where the corresponding wave has a peak.
The {\it decorated permutation} indexing the cell containing $A$ determines the
asymptotic behavior of the soliton solution at $y \to \pm \infty$.
When $t$ is sufficiently small, we can predict the combinatorial structure
of this tropical curve using the {\it $\Le$-diagram} 
indexing the cell containing $A$.
When $A$ comes from a totally positive Schubert cell, 
we show that generically this tropical curve is a realization
of one of Postnikov's  {\it reduced plabic graphs}.
Furthermore, if we label each region of the complement of the tropical curve 
with the dominant exponential in the $\tau$-function, then the labels of the unbounded regions 
form the {\it Grassmann necklace} indexing the cell containing $A$.
Finally, when $A$ belongs to the totally positive  Grassmannian,
we show that the dominant exponentials labeling regions of the tropical curve 
form a {\it cluster} for the cluster algebra of the Grassmannian.  
Letting $t$ vary, one may observe {\it cluster transformations}.  

These previously undescribed connections between KP solitons, cluster algebras, and total positivity
promise to be very powerful.  
For example, using some machinery from total positivity and cluster algebras, 
we solve the {\it inverse problem} for soliton solutions from the totally positive
Grassmannian.

\section{Total positivity for the Grassmannian}

The real Grassmannian $Gr_{kn}$ is the space of all
$k$-dimensional subspaces of $\R^n$.  An element of
$Gr_{kn}$ can be represented by a full-rank $k\times n$ matrix modulo left
multiplication by nonsingular $k\times k$ matrices.  

Let $\binom{[n]}{k}$ be the set of $k$-element subsets of $[n]:=\{1,\dots,n\}$.
For $I\in \binom{[n]}{k}$, let $\Delta_I(A)$
denote the maximal minor of a $k\times n$ matrix $A$ located in the column set $I$.
The map $A\mapsto (\Delta_I(A))$, where $I$ ranges over $\binom{[n]}{k}$,
induces the {\it Pl\"ucker embedding\/} $Gr_{kn}\hookrightarrow \mathbb{RP}^{\binom{n}{k}-1}$,
and the $\Delta_I(A)$ are called {\it Pl\"ucker coordinates}.

\begin{definition} 
The \emph{totally non-negative Grassmannian} $(Gr_{kn})_{\geq 0}$
(respectively, \emph{totally positive Grassmannian} $\Grkntop$)
is the subset of  $Gr_{kn}$
that can be represented by $k\times n$ matrices $A$
with all 
$\Delta_I(A)$ non-negative (respectively, positive).
\end{definition}

Postnikov \cite{Postnikov} gave a decomposition of $(Gr_{kn})_{\geq 0}$ into {\it positroid 
cells.}
For $\M\subseteq \binom{[n]}{k}$,
the {\it positroid cell\/} $S_\mathcal{M}^{tnn}$ is
the set of elements of $(Gr_{kn})_{\geq 0}$ represented by all $k\times n$ matrices $A$ with
the $\Delta_I(A)>0$
for $I\in \mathcal{M}$ and 
$\Delta_J(A)=0$, for $J\not\in \mathcal{M}$.

Clearly $(Gr_{kn})_{\geq 0}$ is a disjoint union of the positroid cells 
$S_\mathcal{M}^{tnn}$ -- in fact it is a 
CW complex \cite{PSW}. Note that
$\Grkntop$ is a positroid cell; it is the unique 
positroid cell in $\Grkn$ of top dimension $k(n-k)$.
Postnikov showed that the cells of $(Gr_{kn})_{\geq 0}$ are naturally labeled 
by (and in bijection
with) the following combinatorial objects \cite{Postnikov}:
\begin{itemize}
\item Grassmann necklaces $\mathcal I$ of type $(k,n)$
\item decorated permutations $\pi$ on $n$ letters with $k$ weak excedances 
\item equivalence classes of {\it reduced plabic graphs} of type $(k,n)$
\item $\Le$-diagrams of type $(k,n)$.
\end{itemize}

For the purpose of studying solitons, we are  interested only 
in the subset of positroid cells which
are {\it irreducible}.

\begin{definition}
We say that a positroid cell 
$S_{\mathcal M}$ is \emph{irreducible} if 
the reduced-row echelon matrix $A$ of any point in the cell
has the following properties:
\begin{enumerate}
\item Each column of $A$ contains at least one nonzero element.
\item Each row of $A$ contains at least one nonzero element in 
addition to the pivot.
\end{enumerate}
\end{definition}

The irreducible positroid cells 
are indexed by:
\begin{itemize}
\item irreducible Grassmann necklaces $\mathcal I$ of type $(k,n)$
\item derangements $\pi$ on $n$ letters with $k$  excedances 
\item equivalence classes of {\it irreducible 
reduced plabic graphs} of type $(k,n)$
\item irreducible $\Le$-diagrams of type $(k,n)$.
\end{itemize}

We now review the definitions of  these objects and 
some of the bijections among them.

\begin{definition}
An  \emph{irreducible Grassmann necklace of type 
$(k,n)$} is a sequence $\mathcal I = (I_1,\dots,I_n)$ of 
subsets $I_r$ of $[n]$ of size $k$  such that, for $i\in [n]$, 
$I_{i+1}=(I_i \setminus \{i\}) \cup \{j\}$ for some $j\neq i$.
(Here indices $i$
are taken modulo $n$.)  
\end{definition}

\begin{example}\label{ex1}
An example of a Grassmann necklace of type $(4,9)$ is 
$(1257, 2357, 3457, 4567, 5678, 6789, 1789, 1289, 1259)$.
\end{example}

\begin{definition}
A \emph{derangement} $\pi=(\pi_1,\dots,\pi_n)$
is a permutation $\pi \in S_n$ which has no fixed points.
An {\it excedance} of $\pi$ 
is a pair $(i,\pi_i)$ such that 
$\pi_i>i$.  We call $i$ the {\it excedance position} and 
$\pi_i$ the {\it excedance value}.  Similarly,
a {\it nonexcedance} is a pair $(i,\pi_i)$ such that $\pi_i<i$.
\end{definition}

\begin{definition}
A \emph{plabic graph}
is a planar undirected graph $G$ drawn inside a disk
with $n$ \emph{boundary vertices\/} $1,\dots,n$ placed in counterclockwise
order
around the boundary of the disk, such that each boundary vertex $i$ 
is incident
to a single edge.\footnote{The convention of \cite{Postnikov} 
was to place the boundary vertices in clockwise order.}  Each internal vertex
is colored black or white.  
\end{definition}

\begin{figure}[h]
\begin{center}
\includegraphics[height=1.05in]{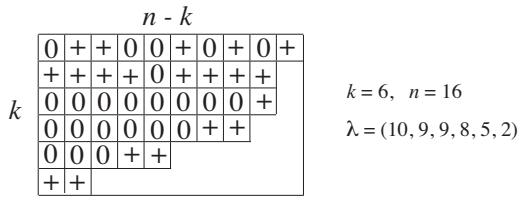}
\par
\end{center}
\caption{A Le-diagram $(\lambda,D)_{k,n}$.\label{LeDiagram}}
\par
\end{figure}

\begin{definition}
Let 
$Y_{\lambda}$ denote the Young diagram of the partition $\lambda$.  A {\it $\Le$-diagram}
(or Le-diagram) 
$L=(\lambda, D)_{k,n}$ of type $(k,n)$ 
is a Young diagram $Y_{\lambda}$ contained in a $k \times (n-k)$ rectangle
together with a filling $D: Y_{\lambda} \to \{0,+\}$ which has the
{\it $\Le$-property}:
there is no $0$ which has a $+$ above it in the same column and a $+$ to its
left in the same row.
A $\Le$-diagram is \emph{irreducible} if each row and each column 
contains at least one $+$.
\end{definition}
See  Figure \ref{LeDiagram} for
an example of an irreducible  $\Le$-diagram.
\begin{theorem}\cite[Theorem 17.2]{Postnikov}\label{necklace}
Let $S_{\mathcal M}^{tnn}$ be a positroid
cell in $(Gr_{kn})_{\geq 0}$. 
For $1 \leq r \leq n$, let $I_r$ be the index set 
of the minor in $\mathcal M$ which is lexicographically
minimal with respect to the order
$r < r+1 <  \dots < n < 1 < 2 < \dots r-1$.
Then $\mathcal I(\mathcal M):=(I_1,\dots,I_n)$ is a Grassmann necklace
of type $(k,n)$.
\end{theorem}

\begin{lemma}\cite[Lemma 16.2]{Postnikov}\label{Postnikov-permutation}
Given an irreducible Grassmann necklace $\mathcal I$, define
a derangement  $\pi=\pi(\mathcal I)$ by requiring that: 
if $I_{i+1} = (I_i \setminus \{i\}) \cup \{j\}$
for $j \neq i$, then $\pi(j)=i$. 
\footnote{Actually Postnikov's convention was to set $\pi(i)=j$ above,
so the permutation we are associating is the inverse one to his.}
Indices are taken modulo $n$.
Then $\mathcal I \to \pi(\mathcal I)$
is a bijection from irreducible
Grassmann necklaces $\mathcal I=(I_1,\dots,I_n)$ of type $(k,n)$
to derangements $\pi(\mathcal I) \in S_n$ with $k$ excedances.
The excedances of 
$\pi(\mathcal I)$ are in positions $I_1$.
\end{lemma}

\begin{remark}
If the positroid cell $S_{\mathcal M}^{tnn}$ is indexed by the Grassmann
necklace $\mathcal I$, the derangement $\pi$, and 
the $\Le$-diagram $L$, then we also refer to this cell as
$S_{\mathcal I}^{tnn}, S_{\pi}^{tnn}$ and $S_{L}^{tnn}$.  The bijections
above preserve the indexing of cells, that is, 
$S_{\mathcal M}^{tnn} = S_{{\mathcal I}({\mathcal M})}^{tnn} 
 = S_{\pi({\mathcal I}({\mathcal M}))}^{tnn}$. 
\end{remark}

\section{Soliton solutions to the KP equation}\label{soliton-background}

Here we explain how to obtain a soliton solution to the KP equation from a point of 
$\Grkn$.

\subsection{From the Grassmannian to the $\tau$-function}

We start by fixing real parameters $\kappa_j$ such that 
$\kappa_1~<~\kappa_2~\cdots~<\kappa_n,$
which  are generic, in the sense
that 
the sums $\sum_{m=1} ^d\kappa_{j_m}$ are all distinct for $2\le d\le k$.

Let $\{E_j;j=1,\ldots,n\}$ be a set of exponential functions in $(x,y,t)\in \mathbb{R}^3$ defined by
\[
E_j(x,y,t):=\exp\left(\kappa_jx+\kappa_j^2y+\kappa_j^3t\right).
\]
If $E_i^{(j)}$ denotes $\partial^jE_i/\partial x^j=\kappa_i^jE_i$, then
the Wronskian determinant with
respect to $x$  of $E_1,\dots,E_n$ is defined by 
\[
{\rm Wr}(E_1,\ldots,E_n)=\det [(E_i^{(j-1)})_{1\le i,j\le n}]
=\prod_{i<j}(\kappa_j-\kappa_i)\,E_1\cdots E_n.
\]

Let $A$ be a full rank $k \times n$ matrix.  We define 
a set of functions $\{f_1,\ldots,f_k\}$ by 
\[
(f_1,f_2,\ldots,f_k)^T = A\cdot (E_1,E_2,\ldots, E_n)^T,
\]
where $(\ldots)^T$ denotes the transpose of the vector $(\ldots)$.
The 
{\it $\tau$-function} of $A$ is defined by
\begin{equation}\label{tauA}
\tau_A(x,y,t):={\rm Wr}(f_1,f_2,\ldots,f_k).
\end{equation}
It's easy to verify that $\tau_A$ only depends on which point of 
$\Grkn$ the matrix $A$ represents.

Applying the Binet-Cauchy identity to the fact that  $f_i=\sum_{j=1}^na_{ij}E_j$ for $i=1,\ldots,k$,
we get \begin{equation}\label{tau}
\tau_A(x,y,t)=\sum_{I\in\binom{[n]}{k}}\Delta_I(A)\,E_I(x,y,t),
\end{equation}
where $E_I(x,y,t)$ with $I=\{\kappa_{j_1},\ldots,\kappa_{j_k}\}$ is defined by
\[
E_I:={\rm Wr}(E_{j_1},E_{j_2},\ldots, E_{j_k})=\prod_{l<m}(\kappa_{j_m}-\kappa_{j_l})\,E_{j_1}\cdots E_{j_k}\,>0.
\]
Therefore if  $A\in \Grkn$, then 
$\tau_A>0$ for all $(x,y,t)\in\mathbb{R}^3$.  

Thinking of $\tau_A$ as a function of $A$, we note 
from \eqref{tau} that
the $\tau$ function encodes the information of the 
Pl\"ucker
embedding.  More specifically, 
if we identify each function $E_I$
with $I=\{j_1,\ldots,j_k\}$  
with the wedge product $E_{j_1}\wedge\cdots\wedge E_{j_k}$, 
then the map  $\tau:Gr_{kn}\hookrightarrow \mathbb{RP}^{\binom{n}{k}-1}$,
$A\mapsto \tau_A$ has the Pl\"ucker coordinates as coefficients.

\subsection{From the $\tau$-function to solutions of the KP equation}

The KP equation
\[
\frac{\partial}{\partial x}\left(-4\frac{\partial u}{\partial t}+6u\frac{\partial u}{\partial x}+\frac{\partial^3u}{\partial x^3}\right)+3\frac{\partial^2 u}{\partial y^2}=0
\]
was proposed by Kadomtsev and Petviashvili in 1970 \cite{KP70}, in order to 
study the stability of the one-soliton solution of the Korteweg-de Vries (KdV) equation
under the influence of weak transverse perturbations.  The KP equation also gives an excellent model
to describe shallow water waves \cite{K10}.

It is well known (see e.g. \cite{H04}) that  the $\tau$-function defined in \eqref{tauA} provides a soliton solution of the KP equation,
\begin{equation}\label{KPsolution}
u_A(x,y,t)=2\frac{\partial^2}{\partial x^2}\ln\tau_A(x,y,t).
\end{equation}
Note that if $A\in \Grkn$, then $u_A(x,y,t)$ is regular.  

\section{From soliton solutions to soliton graphs}
\begin{figure}[h]
\begin{center}
\includegraphics[height=1.1in]{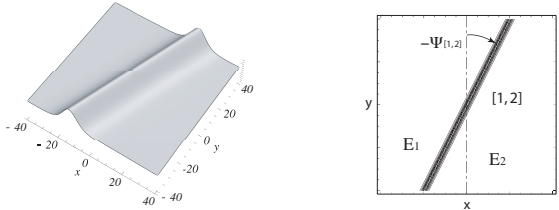}
\par
\end{center}
\caption{A line-soliton solution from $A=(1,1) \in (Gr_{1,2})_{\geq 0}$.
The left figure shows the 3-dimensional profile
of $u_A(x,y,0)$; the right one shows level sets of $u_A(x,y,0)$.
$E_i$ represents the dominant exponential in the region.
\label{fig:1soliton}}
\end{figure}

One can visualize such a solution $u_A(x,y,t)$ in the $xy$-plane 
by drawing level sets of the solution for each time $t$.  
For each $r\in \mathbb{R}$, we denote the corresponding level set by 
\[
C_r(t):=\{(x,y)\in\mathbb{R}^2: u_A(x,y,t)=r\}. 
\]
Figure \ref{fig:1soliton} depicts both a three-dimensional image of a solution $u_A(x,y,t)$,
as well as multiple level sets $C_r(0)$. 
Note that these levels sets are lines parallel to the line of the wave peak.

To study the behavior of $u_A(x,y,t)$ for $A\in S_{\mathcal{M}}^{tnn}$,  we set 
\begin{align*}
\hat{f}_A(x,y,t)&= 
\underset{J\in\mathcal{M}}\max\{\Delta_J(A) E_J(x,y,t)\}\\
              &=\underset{J\in\mathcal{M}}\max \biggl\{e^{\ln\bigl(\Delta_J(A) K_J \bigr) +
                     \sum_{i=1}^k (\kappa_{j_i} x + \kappa^2_{j_i} y +\kappa^3_{j_i} t)}\biggr\},
\end{align*}
where $K_J:=\prod_{\ell < m} (\kappa_{j_m}-\kappa_{j_{\ell}})>0$.
From \eqref{tau}, we see that generically, $\tau_A$ can be approximated by $\hat{f}_A$.

Let $f_A(x,y,t)$ be the closely related function 
\begin{equation}
{f}_A(x,y,t) =
              \underset{J\in\mathcal{M}}\max \left\{\ln\bigl(\Delta_J(A)K_J\bigr) +
                     \sum_{i=1}^k (\kappa_{j_i} x + \kappa^2_{j_i} y +\kappa^3_{j_i} t)\right\}.
\end{equation}
Clearly a given term dominates $f_A(x,y,t)$ if and only if its exponentiated version dominates
$\hat{f}_A(x,y,t)$.

\begin{definition}\label{contour}
Given a solution $u_A(x,y,t)$ of the KP equation as in 
\eqref{KPsolution}, we define its \emph{contour plot}  $\CC_{t_0}(u_A)$
for each $t=t_0$ to be the locus in $\R^2$ where $f_A(x,y,t=t_0)$ is not 
linear. 
\end{definition}

\begin{remark}
$\CC_{t_0}(u_A)$ provides an approximation of the location of the wave crests.
\end{remark}

It follows from Definition \ref{contour} that $\CC_{t_0}(u_A)$
is a one-dimensional piecewise linear subset of the $xy$-plane.

\begin{proposition}
If each $\kappa_{i}$ is an integer, then 
 $\CC_{t_0}(u_A)$ is a tropical curve
in $\R^2$.  
\end{proposition}

Note that each region of the complement of $\CC_{t_0}(u_A)$ in 
$\R^2$ is a domain of linearity for $f_A(x,y,t_0)$, 
and hence each region is  
naturally associated to a {\it dominant exponential}  
$\Delta_J(A) E_J(x,y,t_0)$ from the $\tau$-function \eqref{tau}.
We call the line segments comprising 
$\CC_{t_0}(u_A)$ {\it line-solitons}.\footnote{In general, there exist phase-shifts
which also appear as line segments (see \cite{CK3}). However the phase-shifts
depend only on the $\kappa$ parameters, and we ignore them in this paper.}  Some of 
these line-solitons have finite length, while others are unbounded
and extend in the $y$ direction to $\pm \infty$.  We call these
{\it unbounded} line-solitons.
Note that each line-soliton
represents a balance between two dominant exponentials
in the $\tau$-function.

\begin{lemma}\label{separating}\cite[Proposition 5]{CK3}
The dominant exponentials of the $\tau$-function in adjacent regions
of the contour plot in the $xy$-plane are of the form $E(i,m_2,\dots,m_k)$ and 
$E(j, m_2,\dots,m_k)$.  
\end{lemma}

According to Lemma \ref{separating}, 
those two exponential terms have $k-1$ common
phases, 
so we call the soliton separating them a {\it line-soliton of type $[i,j]$}.
Locally we have 
\begin{align*}
\tau_A&\approx \Delta_I(A)E_I+\Delta_J(A)E_J\\
&=\left(\Delta_I(A)K_I E_i+\Delta_J(A)K_J E_j\right)\prod_{l=2}^kE_{m_l}
\end{align*}
with $K_I=\prod_{j=2}^k|\kappa_i-\kappa_{m_j}|\prod_{l<j}|\kappa_{m_j}-\kappa_{m_l}|$,
so the
equation for this line-soliton is
\begin{equation}\label{eq-soliton}
x+(\kappa_i+\kappa_j)y+(\kappa_i^2+\kappa_i\kappa_j+\kappa_j^2)t=\frac{1}{\kappa_j-\kappa_i}\ln\frac{\Delta_I(A)K_I}{\Delta_J(A)K_J}.
\end{equation}
Note that the ratio of the Pl\"ucker coordinates labeling the regions
separated by the line-soliton determines the location of the line-soliton.

\begin{remark}\label{slope}
Consider a line-soliton given by  \eqref{eq-soliton}.
Compute the angle $\Psi_{[i,j]}$ 
between the line-soliton and the positive $y$-axis,
measured in the counterclockwise direction, so that the negative $x$-axis
has an angle of $\frac{\pi}{2}$ and the positive $x$-axis has an 
angle of $-\frac{\pi}{2}$. Then $\tan \Psi_{[i,j]} = \kappa_i+\kappa_j$.
Therefore we refer to $\kappa_i+\kappa_j$ as the \emph{slope} of the 
$[i,j]$ line-soliton (see Figure \ref{fig:1soliton}).
\end{remark}

We will be interested in the combinatorial structure of a contour plot,
that is, the pattern of how line-solitons interact with each other.  To 
this end, in Definition \ref{soliton-graph} 
we will associate a {\it soliton graph} to each contour plot.
 
Generically we expect a point of a contour plot at which several line-solitons
meet to have degree $3$; we regard such a point as a trivalent 
vertex.   Three line-solitons meeting 
at a trivalent
vertex exhibit a {\it resonant interaction} (this corresponds
to the {\it balancing condition} for a tropical curve).
One may also have two line-solitons which cross over
each other, forming an $X$-shape: we call this
an \emph{$X$-crossing}, but do not regard it as a vertex.
In general, there exists a phase-shift at each $X$-crossing. However
 we ignore them in this paper as explained in the footnote 3.
Vertices of degree greater than $4$
are also possible.  

\begin{definition}\label{def:generic}
A contour plot is called \emph{generic} if all interactions of line-solitons
are at trivalent vertices or are $X$-crossings.
\end{definition}

The following definition of {\it soliton graph} 
forgets the metric 
data of the contour plot, but preserves
the data of how line-solitons interact and which exponentials are dominant.

\begin{definition}\label{soliton-graph}
Let $\CC_{t_0}(u_A)$ be a generic contour plot 
with $n$ unbounded
line-solitons.
Color a trivalent 
vertex black (respectively, white) 
if it has a unique edge extending downwards (respectively, upwards) from it.
Label each region with the dominant exponential $E_I$
and each edge (line-soliton) by the {\it type} $[i,j]$ of that
line-soliton.
Preserve the topology of the metric graph, 
but forget the metric structure. 
Embed the resulting graph with bicolored vertices and $X$-crossings 
into a disk with $n$ boundary vertices, 
replacing each unbounded line-soliton with an edge 
that ends at a boundary vertex.
We call this labeled graph
a \emph{soliton graph}.
\end{definition}

\begin{figure}[h]
\begin{center}
\includegraphics[height=1.8in]{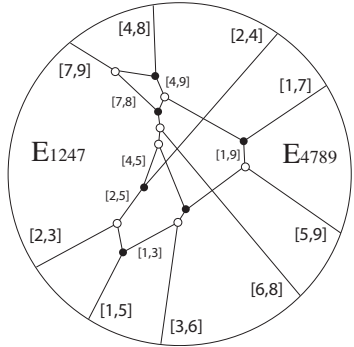}
\end{center}
\caption{Example of a soliton graph associated to
$\S_{\pi}^{tnn}$ with
$\pi=(7,4,2,8,1,3,9,6,5)$.  Each $E_{ijkl}$ represents
the dominant exponential in the $\tau$-function.
This soliton graph was obtained from a contour plot by 
embedding it in a disk and coloring
vertices appropriately. 
\label{contour-soliton}}
\end{figure}

See Figure \ref{contour-soliton} for an example of 
a soliton graph.
Although we have not labeled all regions or all edges, 
the remaining labels can be determined using Lemma \ref{separating}.

\section{Permutations and soliton asymptotics}\label{permutations}
Given a contour plot $\CC_{t_0}(u_A)$ where $A$ belongs
to an irreducible positroid cell and $t_0$ is arbitrary, we show that the 
labels of the unbounded solitons 
allow us to determine which positroid cell $A$ belongs to.
Conversely, given $A$ in the irreducible 
positroid cell $S_{\pi}^{tnn}$, 
we can predict the asymptotic behavior of the unbounded solitons
in $\CC_{t_0}(u_A)$.
\begin{theorem}
Suppose $A$ is an element of an irreducible positroid cell in 
$(Gr_{kn})_{\geq 0}$.
Consider the contour plot $\CC_{t_0}(u_A)$ for any time $t_0$.
Then there are $k$ unbounded
line-solitons at $y\gg0$ which are labeled by 
pairs $[e_r,j_r]$ with $e_r<j_r$, and there are $n-k$ 
unbounded line-solitons at $y\ll0$ which are labeled by pairs
$[i_r,g_r]$ with $i_r<g_r$.  We obtain a derangement in $S_n$
with $k$ excedances by setting $\pi(e_r) = j_r$ and 
$\pi(g_r) = i_r$.  Moreover, $A$ must be an element of 
the cell $S_{\pi}^{tnn}$.
\end{theorem}

The first part of this theorem follows from work of 
Chakravarty and Kodama
\cite[Prop. 2.6 and 2.9]{CK1}, \cite[Theorem 5]{CK3}.
Our contribution  is that the derangement $\pi$ is precisely 
the derangement labeling the cell $S_{\pi}^{tnn}$ that 
$A$ belongs to.
This fact is the first step towards establishing that various other
combinatorial objects in bijection with positroid cells
(Grassmann necklaces, plabic graphs) carry useful information about the 
corresponding soliton solutions.

We now give a concrete algorithm for writing down the asymptotics 
of the soliton solutions of the KP equation.
\begin{figure}
\begin{center}
\includegraphics[height=3.6cm]{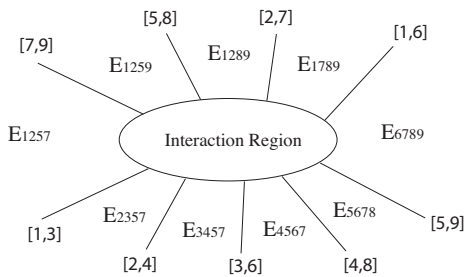}
\end{center}
\caption{Unbounded line-solitons for $\pi=(6,7,1,2,8,3,9,4,5)$.
Each $E_{ijkl}$ shows the dominant exponential in this region.
\label{671283945}}
\end{figure}
\begin{theorem}\label{algo}
Fix generic parameters $\kappa_1<\dots < \kappa_n$.  
Let $A$ be an element from an irreducible 
positroid cell
$\SSS_{\pi}^{tnn}$ in $(Gr_{kn})_{\geq 0}$.
(So $\pi$ must have $k$ excedances.)
For any $t_0$, the asymptotic behavior of the contour plot $\mathcal{C}_{t_0}(u_A)$
--  i.e. its unbounded line-solitons, 
and the dominant exponentials in its unbounded regions -- can  be
read off from $\pi$ as follows.
\begin{itemize}
\item For $y \gg0$, there is an unbounded line-soliton of type
$[i,\pi(i)]$ for each
excedance $\pi(i)>i$.  
From left to right,
list these solitons in decreasing order of the quantity
$\kappa_i + \kappa_{\pi(i)}$.
\item For $y\ll 0$, there is an unbounded line-solitons  of type
$[\pi(j),j]$ for each nonexcedance $\pi(j)<j$.  From left to right,
list these solitons in increasing order of 
$\kappa_j + \kappa_{\pi(j)}$.
\item Label the unbounded region for $x\ll 0$ with 
the exponential $E_{i_1,\dots,i_k}$,
where $i_1,\dots,i_k$ are the excedance positions of $\pi$.
\item Use Lemma \ref{separating} to label the remaining
unbounded regions of the contour plot.
\end{itemize}
\end{theorem}

\begin{example}\label{ex4}
Consider the positroid cell corresponding to 
$\pi = (6,7,1,2,8,3,9,4,5)\in S_9$.  The algorithm of Theorem 
\ref{algo} gives rise to the picture in Figure \ref{671283945}.
If one reads the dominant exponentials in 
counterclockwise order,
starting from the region at the left, then one recovers 
the Grassmann necklace $\mathcal I$ from Example \ref{ex1}.
Also note that $\pi(\mathcal I) = \pi$.
See Theorem \ref{necklace-soliton}.
\end{example}

\section{Grassmann necklaces and  soliton asymptotics}

One particularly nice
class of positroid cells is the {\it TP} or 
{\it totally positive Schubert cells}.
These are the positroid cells indexed by $\Le$-diagrams which 
are filled with all $+$'s, or equivalently, the positroid
cells indexed by derangements $\pi$ such that $\pi^{-1}$ has
at most one descent.
When $S_{\pi}^{tnn}$ is a TP  Schubert cell,
we can make a link between the corresponding soliton solutions of the 
KP equation and Grassmann necklaces.

\begin{theorem}\label{necklace-soliton}
Let $A$ be an element of 
a TP Schubert 
cell $\mathcal{S}_{\pi}^{tnn}$, and consider the contour
plot $\CC_{t_0}(u_A)$ for an arbitrary time $t_0$.
Let the index sets of the dominant exponentials of 
the unbounded regions of $\CC_{t_0}(u_A)$ be 
denoted $R_1, \dots, R_n$,
where $R_1$ labels the region at $x\ll 0$, and $R_2,\dots,R_n$
label the regions in the counterclockwise direction from $R_1$.
Then $(R_1,\dots, R_n)$ is a Grassmann necklace $\mathcal I$
and $\pi(\mathcal I) = \pi$.
\end{theorem}

Theorem \ref{necklace-soliton} is illustrated in Example \ref{ex4}.  

\begin{remark}
Theorem \ref{necklace-soliton} does not hold if we replace ``TP
Schubert cell" by ``positroid cell."  
\end{remark}

\section{From soliton graphs to generalized plabic graphs}

In this section we associate a {\it generalized plabic graph} $Pl(C)$ to 
each soliton graph $C$.  
We then show that from $Pl(C)$ --
whose only labels are on the boundary vertices -- 
we can recover the labels of the line-solitons and dominant exponentials
of $C$.

\begin{definition}
A \emph{generalized plabic graph\/}
is a connected graph embedded in a disk
with $n$ \emph{boundary vertices\/} labeled $1,\dots,n$ placed 
in \emph{any} order
around the boundary of the disk, such that each boundary vertex $i$ 
is incident
to a single edge. Each internal vertex must have degree at least two, and
is colored black or white.  Edges
are allowed to form 
$X$-crossings  (this is not considered to be a vertex). 
\end{definition}

We now generalize the notion
of trip from \cite[Section 13]{Postnikov}.

\begin{definition}\label{gen:trip}
Given a generalized plabic graph $G$,
the \emph{trip} $T_i$ is the directed path which starts at the boundary vertex 
$i$, and follows the ``rules of the road": it turns right at a 
black vertex,  left at a white vertex, and goes straight  through an
$X$-crossing. Note that $T_i$  will also 
end at a boundary vertex.  The \emph{trip permutation} $\pi_G$
is the permutation such that $\pi_G(i)=j$ whenever the trip
starting at $i$ ends at $j$. 
\end{definition}

We use these trips to associate a canonical labeling of edges and regions
to each generalized plabic graph.

\begin{definition}\label{labels}
Given a generalized plabic graph $G$ with $n$ boundary vertices, 
start at each boundary
vertex $i$ and label every edge along trip $T_i$ with $i$.
Such a trip divides the disk containing $G$ into two parts: 
the part to the left of $T_i$, and the part to the right.
Place an $i$ in every region which is to the left of $T_i$.
After repeating this procedure for each boundary vertex,
each edge will be labeled by up to two numbers (between $1$ and $n$),
and each region will be labeled by a collection of numbers.
Two regions separated by an edge labeled $ij$ will have 
region labels $S$ and $(S\setminus \{i\}) \cup \{j\}$.
When an edge is assigned two numbers $i<j$, 
we write $[i,j]$
on that edge, or $\{i,j\}$ or $\{j,i\}$ if we do not wish 
to specify the order of $i$ and $j$.  
\end{definition}

\begin{figure}[h]
\begin{center}
\includegraphics[height=1.9in]{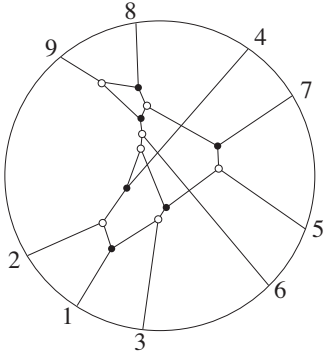}
\end{center}
\caption{Example of a generalized plabic graph $G(C)$.
\label{contour-soliton2}}
\end{figure}

 
\begin{definition}\label{soliton2plabic}
Fix an irreducible cell $\S_{\pi}^{tnn}$ of $(Gr_{kn})_{\geq 0}$.
To each soliton graph $C$ coming from a point of that cell we associate 
a generalized plabic graph
$Pl(C)$ by:
\begin{itemize}
\item labeling the boundary vertex 
incident to the edge $\{i,\pi_i\}$ 
by $\pi_i=\pi(i)$,
\item forgetting the labels of all edges
and regions.
\end{itemize}
\end{definition}
See Figure \ref{contour-soliton2} for 
the generalized plabic graph $Pl(C)$ corresponding to the 
soliton graph $C$ from 
Figure \ref{contour-soliton}.

\begin{theorem}\label{soliton-plabic}
Fix an irreducible cell $\S_{\pi}^{tnn}$ of $(Gr_{kn})_{\geq 0}$,
and consider 
a soliton graph $C$ coming from a point of that cell.
Then the trip permutation associated to the plabic graph $Pl(C)$ is $\pi$,
and by labeling edges and regions of $Pl(C)$ according
to Definition \ref{labels}, we will recover the original
labels in  $C$.
\end{theorem}
We invite the reader to apply Definition \ref{labels} to  
Figure \ref{contour-soliton2},
and then compare the result to Figure \ref{contour-soliton}.

\begin{remark}
By Theorem \ref{soliton-plabic}, we can identify each soliton graph $C$
with its generalized plabic graph $Pl(C)$.
\end{remark}

\section{Soliton graphs for positroid cells when $t\ll 0$}
\label{plabic-soliton}
In this section we give an algorithm for producing a 
generalized plabic graph $G_{-}(L)$ from the $\Le$-diagram $L$ of a positroid cell
$S_{L}^{tnn}$.
It turns out that
this generalized plabic graph gives rise to
the soliton graph for a generic point of 
the cell $\S_{\Le}^{tnn}$, at time $t\ll0$ sufficiently small. 
\begin{algorithm}  \label{LeToPlabic} 
Given  a $\Le$-diagram $L$, construct $G_{-}(L)$ as follows:
\begin{enumerate}
\item Start with a $\Le$-diagram $L$ contained in a $k\times (n-k)$
rectangle. Label its southeast border by the numbers $1$ to $n$, starting from the northeast corner.  Replace $0$'s and $+$'s by ``crosses'' and ``elbows''.  From each label $i$ on the southeast border,
follow the associated ``pipe'' northwest, and label its destination by $i$ as well.
\item Add an edge, and one white and one black vertex  
to each elbow, as shown
in the upper right of Figure \ref{LePlabic}.  Forget the labels
of the southeast border.  If there is an endpoint of a pipe on the east or south border whose pipe
starts by going straight, then erase the straight portion preceding the first elbow.
\item Forget any degree $2$ vertices, and forget
any edges of the graph which end 
at the southeast border of the diagram.
Denote the resulting
graph  $G_-(L)$.  
\item After embedding the graph in a disk 
with $n$ boundary vertices,
we obtain a generalized plabic graph,
which we also denote $G_-(L)$.
If desired, stretch and rotate $G_-(L)$ so that the boundary vertices
at the west side of the diagram are at the north instead.
\end{enumerate}
\end{algorithm}
\begin{figure}[h]
\centering
\includegraphics[height=2in]{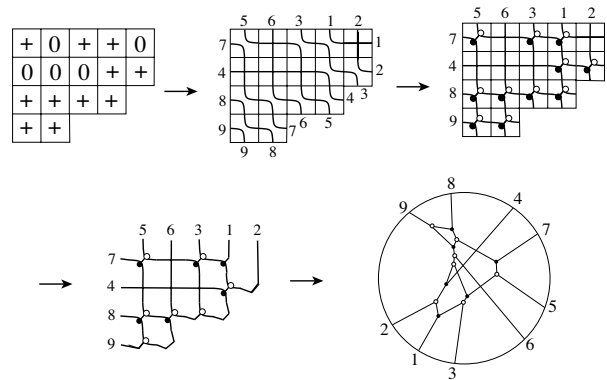}
\caption{Algorithm \ref{LeToPlabic} for $\S_{L}^{tnn}$ with 
$L$ the Le-diagram from the upper left.
\label{LePlabic}}
\end{figure}

Figure \ref{LePlabic} illustrates the steps of Algorithm \ref{LeToPlabic}. 
Note that this produces the graph from Figure 
\ref{contour-soliton2}.


\begin{theorem}\label{t<<0}
Let $L$ be a $\Le$-diagram and $\pi=\pi(L)$. Then  $G_{-}(L)$ 
has
trip permutation $\pi$.  Label its edges and regions according to the rules of the road. When $S_L^{tnn}$ is a TP
Schubert cell, then $G_{-}(L)$ coincides with the soliton graph $G_t(u_A)$, provided that $A\in S_{L}^{tnn}$ is
generic and $t\ll0$ sufficiently small. When $S_L^{tnn}$ is an arbitrary positroid cell,  we can realize
$G_{-}(L)$ as ``most" of a soliton graph $G_t(u_A)$ for $A\in S_{L}^{tnn}$ and $t\ll0$. Moreover, we
can construct $G_t(u_A)$ from $G_{-}(L)$ by extending the unbounded edges of $G_{-}(L)$ and introducing
$X$-crossings as necessary so as to satisfy the conditions of Theorem \ref{algo}.
\end{theorem}

\section{Reduced plabic graphs  and cluster algebras}\label{Reduced-Cluster}

The most important plabic graphs are those which are {\it reduced}
\cite[Section 12]{Postnikov}.  Although it is not easy to 
characterize reduced plabic graphs (they are defined to be 
plabic graphs whose {\it move-equivalence class} contains no 
graph to which one can apply a {\it reduction}), they are 
important because of their application to cluster algebras
and  parameterizations of 
cells.

\begin{theorem}\label{th:reduced}
Let $A$ be a point of a TP Schubert cell, let $t_0$
be an arbitrary time, 
and suppose that the contour plot $\CC_{t_0}(u_A)$ is generic and has 
no X-crossings.
Then the soliton graph associated to $\CC_{t_0}(u_A)$
is a reduced plabic graph.
\end{theorem}

Cluster algebras are a class of commutative rings with a 
remarkable combinatorial structure, which were defined by 
Fomin and Zelevinsky \cite{FZ}.  Scott \cite{Scott} proved that
Grassmannians have a cluster algebra structure.

\begin{theorem} \cite{Scott}\label{cluster-Grassmannian}
The coordinate ring of the 
(affine cone over the) Grassmannian has the structure of a 
cluster algebra.  Moreover, the set of labels of the regions
of any reduced plabic graph for the TP (totally positive) Grassmannian
comprises a \emph{cluster} for this cluster algebra.
\end{theorem}


\begin{remark}
Scott's strategy in \cite{Scott} was to show that certain labelings of
{\it alternating strand diagrams} for the TP
Grassmannian gave rise to clusters.  However, alternating strand 
diagrams are in bijection with reduced plabic graphs \cite{Postnikov}, and 
under this bijection,
Scott's labelings of alternating strand diagrams  correspond to 
the labelings of regions of plabic graphs induced by the various
trips in the plabic graph.
\end{remark}

\begin{corollary}\label{soliton-cluster}
The Pl\"ucker coordinates labeling regions of a
generic soliton graph with no X-crossings for $(Gr_{k,n})_{>0}$
form a cluster for the cluster algebra
of $Gr_{k,n}$.
\end{corollary}

Conjecturally, every positroid cell $\S_{\pi}^{tnn}$
of the totally non-negative
Grassmannian also carries a cluster algebra structure, 
and the Pl\"ucker coordinates labeling the regions of any reduced
plabic graph for $\S_{\pi}^{tnn}$ should be a cluster for that cluster
algebra.  In particular, the TP Schubert cells should carry cluster algebra
structures.  Therefore we conjecture that Corollary 
\ref{soliton-cluster} holds with ``TP Schubert cell" replacing ``TP
Grassmannian."  Finally, there should be a suitable generalization
of Corollary \ref{soliton-cluster} for arbitrary 
positroid cells.

\section{The inverse problem}
\label{inverse}

The {\it inverse problem} for soliton solutions of the KP equation
is the following:
given a time $t$ together with the contour plot of a soliton 
solution, can one reconstruct the point of 
$(Gr_{k,n})_{\geq 0}$ which gave rise to the solution?

\begin{theorem}\label{inverse1}
Fix $\kappa_1<\dots<\kappa_n$ as usual.
Consider a generic contour plot of a soliton solution
coming from a point $A$ of a positroid cell $S_{\pi}^{tnn}$,
for $t\ll0$.  Then from the contour plot together with $t$ we 
can uniquely reconstruct the point $A$.
\end{theorem}

The strategy of the proof is as follows. From the contour
plot together with $t$, we can 
reconstruct the value of each of the dominant
exponentials (Pl\"ucker coordinates) labeling regions of the graph.
We have shown how to 
use the $\Le$-diagram to construct
the soliton graph for a positroid cell when 
$t\ll0$ is sufficiently small, which allows us to identify
what is the set of Pl\"ucker coordinates which label regions of the 
graph. We then show that this collection of Pl\"ucker coordinates
contains a subset of Pl\"ucker coordinates, which Talaska \cite{Talaska}
showed were sufficient for reconstructing
the original point of $\S_{\pi}^{tnn}$. 

Using Theorem \ref{inverse1}, Corollary \ref{soliton-cluster},
 and the cluster algebra structure 
for Grassmannians, we can solve the inverse problem for the TP
Grassmannian for any time $t$. 
\begin{theorem}
Consider a generic contour plot of a soliton solution
coming from a point $A$ of the TP Grassmannian,
at an {\it arbitrary} time $t$.  If the contour plot has no X-crossings, 
then from the contour plot together with $t$ we 
can uniquely reconstruct the point $A$.
\end{theorem}

\section{Triangulations of a polygon and soliton graphs}

We now explain how to use triangulations of an $n$-gon to produce 
all soliton graphs for the TP Grassmannian $(Gr_{2,n})_{>0}$.
\begin{algorithm}  \label{TriangulationA} 
Let $T$ be a triangulation of an $n$-gon $P$, whose $n$ vertices are labeled
by the numbers $1,2,\dots,n$, in counterclockwise order. Therefore
each edge of $P$ and each diagonal of $T$ is specified by a pair of distinct integers
between $1$ and $n$. The following 
procedure yields a labeled graph  $\Psi(T)$.
\begin{enumerate}
\item Put a black vertex in the interior of each triangle in $T$.
\item Put a white vertex at each of the $n$ vertices of $P$
which is incident to a diagonal of $T$; put a black vertex
at the remaining vertices of $P$.
\item Connect each vertex which is 
inside a triangle of $T$ to the three vertices
of that triangle.  
\item Erase the edges of $T$, and contract every pair of  adjacent vertices
which have the same color.  This produces a new graph $G$
with $n$ boundary vertices, in bijection with the vertices of the original 
$n$-gon $P$.
\item Add one unbounded ray to each of the boundary vertices of 
$G$, so as to produce a new (planar) graph  $\Psi(T)$.  Note that 
$\Psi(T)$ divides the plane into regions; the bounded regions correspond to 
the diagonals of $T$, and the unbounded regions correspond to the edges of $P$.
\item Resolve any non-trivalent vertices into trivalent vertices.
\end{enumerate}
\end{algorithm}
\begin{theorem}
The graphs $\Psi(T)$ constructed above {\it are} soliton graphs
for $(Gr_{2,n})_{>0}$. Conversely, any generic soliton graph
with no X-crossings
for $(Gr_{2,n})_{>0}$ comes from Algorithm \ref{TriangulationA}.
\end{theorem}

\begin{figure}[h]
\centering
\includegraphics[height=1.1in]{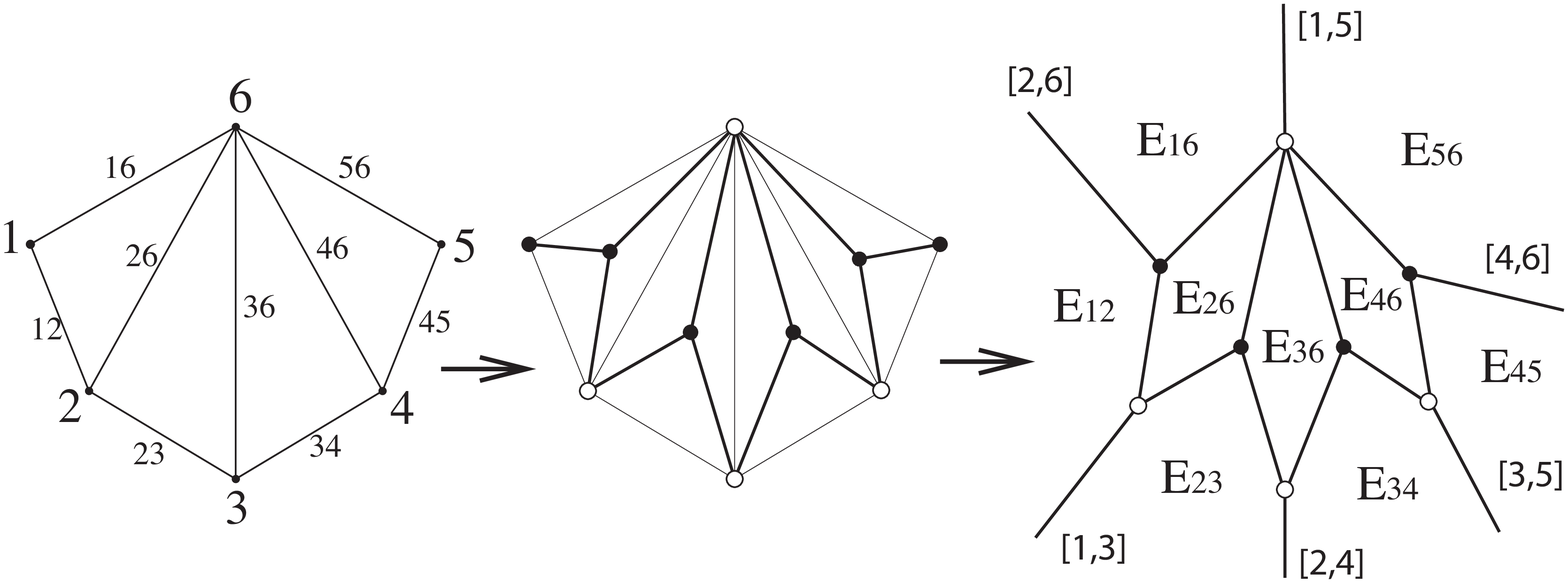}
\caption{Algorithm \ref{TriangulationA}, starting from a triangulation of 
a hexagon.}
\label{Psi}
\end{figure}

Flipping a diagonal in a triangulation corresponds 
to a mutation in the cluster algebra.
In our setting, 
each mutation may be considered as an 
evolution along a 
flow of the KP hierarchy 
defined by the symmetries of
the KP equation.
\begin{acknowledgments}
The first author was partially supported by the NSF grant DMS-0806219,
and the second author was partially supported by the NSF grant
DMS-0854432 and a Sloan fellowship.
\end{acknowledgments}

\end{article}

\begin{thebibliography}{2}
\bibitem{KP70}
    B. B. Kadomtsev and V. I. Petviashvili,
    On the stability of solitary waves in weakly dispersive media,
    {\it Sov. Phys. - Dokl.} {\bf 15} (1970) 539-541.

\bibitem{Sato} M. Sato, Soliton equations as dynamical systems on 
an infinite dimensional Grassmannian manifold, RIMS Kokyuroku (Kyoto 
University) 439 (1981), 30--46.

\bibitem{K04}
    Y. Kodama,
    Young diagrams and $N$-soliton solutions of the KP equation,
    {\it J. Phys. A: Math. Gen.}, {\bf 37} (2004) 11169-11190.
   
\bibitem{BC06}
    G. Biondini and S. Chakravarty, 
    Soliton solutions of the Kadomtsev-Petviashvili II equation,
    {\it J. Math. Phys.}, {\bf 47} (2006) 033514 (26pp).

\bibitem{CK1} S. Chakravarty, Y. Kodama, Classification of the line-solitons of KPII, 
J. Phys. A: Math. Theor. 41 (2008)
275209 (33pp).


\bibitem{CK3} S. Chakravarty, Y. Kodama, Soliton solutions of the KP equation
and applications to shallow water waves, Stud. Appl. Math. 123 (2009)
83--151.

\bibitem{Lusztig2} G. Lusztig, Total positivity in partial flag manifolds,
Represent. Theory 2 (1998), 70-78.

\bibitem{Lusztig3}
G. Lusztig, Total positivity in reductive groups, in:
Lie theory and geometry: in honor of Bertram Kostant,
Progress in Mathematics 123, Birkhauser, 1994.
\bibitem{Postnikov} A. Postnikov, 
Total positivity, Grassmannians, and networks,
arXiv:math.CO/060976v1.
\bibitem{PSW} A. Postnikov, D. Speyer, L. Williams, Matching polytopes,
toric geometry, and the non-negative part of the Grassmannian,
J. Alg. Combin., 30 (2009), 173--191.

\bibitem{K10}
     Y. Kodama,
     KP soliton in shallow water,
     {\it J. Phys. A: Math. Theor.} {\bf 43} (2010) 434004 (54pp).    
    

\bibitem{H04}
    R Hirota,
   \textit{The Direct Method in Soliton Theory}
   (Cambridge University Press, Cambridge, 2004), Chapter 3.
  







  
   
\bibitem{FZ} S. Fomin, A. Zelevinsky, Cluster Algebras I: Foundations,
\emph{J. Amer. Math. Soc.}, 15 (2002), 497--529.    

      
    
     




\bibitem{Scott} J. Scott, Grassmannians and cluster algebras, 
Proc. London Math. Soc. (3) 92 (2006) 345--380.


\bibitem{Talaska} K. Talaska, Combinatorial formulas for $\Le$-coordinates
in a totally nonnegative Grassmannian, 
arXiv:0812.0640, to appear in J. Combin. Theory Ser. A.


     

 

    






\end{thebibliography}
\end{document}